\documentstyle[amstex,amssymb,12pt]{amsart}

\newcommand{\cplx}{{\Bbb C}}
\newcommand{\cz}{\cplx[z_1^{\pm 1},\dots,z_N^{\pm 1}]}
\newcommand{\zint}{{\Bbb Z}}
\newcommand{\zintp}{{\Bbb Z}_{\geq 0}}
\newcommand{\cplxn}{\cplx^{n}}
\newcommand{\ep}{\epsilon}
\newcommand{\Hhat}{\widehat{H}_N(q)}
\newcommand{\HH}{H_N(q)}
\newcommand{\sll}{\widehat{ \frak{s}\frak{l}}_n}
\newcommand{\End}{{\mathrm E}{\mathrm n}{\mathrm d}}
\newcommand{\Ker}{{\mathrm K}{\mathrm e}{\mathrm r}}

\newcommand{\Uslt}{U_q'(\widehat{ \frak{s}\frak{l}}_{2})}
\newcommand{\sln}{ {{\frak s} {\frak l}}_n}

\newcommand{\glN}{ {{\frak g} {\frak l}}_N}
\newcommand{\slt}{\widehat{ \frak{s}\frak{l}}_2}
\newcommand{\tl}{ \frak{s}\frak{l}_2}

\newcommand{\mm}{{\bold m}}

\newcommand{\zz}{{\bold z}}
\newcommand{\ee}{{\bold e}}

\newcommand{\oro}{\overline{\rho}^{M,k}_{l}}
\newcommand{\orom}{\overline{\rho}^{M,k}_{l,m}}

\newcommand{\MC}{{\cal M}}
\newcommand{\EC}{{\cal E}}

\newcommand{\MN}{\tilde {\cal M}_{N}}
\newcommand{\setN}{\{1,2,\dots,N \}}
\newcommand{\Ph}{\Phi_{\s}^{\l} (z) }

\newcommand{\Php}{\tilde{\Phi }_{\s}^{\l} (z) }

\newcommand{\Phm}{\tilde{\Phi }_{min}^{\l}(z) }
\renewcommand{\l}{\lambda}
\newcommand{\s}{\sigma}

\newcommand{\ba}{\begin{array}}
\newcommand{\ea}{\end{array}}
\newcommand{\beq}{\begin{equation}}
\newcommand{\eeq}{\end{equation}}
\newcommand{\bqa}{\begin{eqnarray}}
\newcommand{\eqa}{\end{eqnarray}}
\newcommand{\bqas}{\begin{eqnarray*}}
\newcommand{\eqas}{\end{eqnarray*}}
\newcommand{\halmos}{\rule{5pt}{5pt}}

\numberwithin{equation}{section}
\newtheorem{df}{\bf Definition}
\newtheorem{prop}{\bf Proposition}
\newtheorem{thm}[prop]{\bf Theorem}
\newtheorem{lemma}{\bf Lemma}
\newtheorem{cor}{\bf Corollary}

\newenvironment{rmk}{\noindent{\bf Remark}\hskip 5pt}{}

\renewenvironment{pf}{\noindent {\em {\normalsize P}\normalsize{roof.}}  
\normalsize\hskip 5pt}{\hfill\halmos}

\begin{document}
 
\title[The decomposition of level--$1$ modules]
 {The decomposition of level--$1$ irreducible highest weight modules with respect to the level--$0$ actions of the quantum affine algebra $U'_q(\sll )$}
\author{Kouichi Takemura}
\address[ ]{Research Institute for Mathematical Sciences,
Kyoto University, Japan   }
\email{ takemura@@kurims.kyoto-u.ac.jp}  
\thanks{K.T. is supported by the JSPS Research Fellowship for Young Scientist.}
\maketitle

\begin{abstract}
We decompose the level--$1$ irreducible highest weight modules of the quantum affine algebra $U_q(\sll )$ with respect to the level--$0$ $U'_q (\sll )$--action defined in \cite{STU}. The decomposition is parameterized by the skew Young diagrams of the border strip type.
\end{abstract}

\section{Introduction}

In the papers \cite{VV2} and \cite{STU} it was shown that the $q$-deformed Fock space module of the quantum affine algebra $U_q(\sll )$ admits an action of a new remarkable object -- the so-called quantum toroidal algebra introduced in \cite{GKV} and \cite{VV} as a $q$-deformation of the universal central extension of the $\sln$-valued double-loop Lie algebra. The action of the quantum toroidal algebra on the $q$-Fock space depends on two parameters: the deformation parameter $q$ and an extra parameter $p,$ when  values of these parameters are taken to be generic complex numbers, the $q$-Fock space is known to be irreducible with respect to this action.

The quantum toroidal algebra has two subalgebras, $U'_q (\sll )^{(1)}$ and $U'_q (\sll )^{(2)},$ both isomorphic to $U'_q (\sll ).$ Accordingly, the $q$-Fock space admits two $U'_q (\sll )$-actions. 

The first of these actions has level 1, and coincides with the action originally introduced  by Hayashi in \cite{H}. The irreducible decomposition of the   $q$-Fock space with respect to this action was  given in \cite{KMS} by using the semi-infinite $q$-wedge construction due to \cite{St}. 

The second of the $U'_q (\sll )$-actions has level 0, the irreducible decomposition of the $q$-Fock space with respect to this action was constructed in \cite{TU} at generic values of the parameters $q$ and $p.$

Kashiwara, Miwa and Stern \cite{KMS} have shown, that the level-1 action of $U'_q (\sll )$ on the $q$-Fock space is centralized by the action of the Heisenberg algebra. In the paper \cite{STU} it was proven that the proper ideal of the $q$-Fock space  generated by the negative-frequency part of the Heisenberg algebra is invariant under the action of the  quantum toroidal algebra provided the value of the parameter $p$ in the latter is set to be equal to 1. The quotient of the $q$-Fock space  by this ideal is isomorphic to one of the irreducible level-1 highest weight modules of  $U_q (\sll ).$ As a consequence, each of these modules admits an action of the quantum toroidal algebra.  

The corresponding action of the subalgebra $U'_q (\sll )^{(1)}$ is irreducible, it is just the standard level-1 action on the highest weight irreducible module of $U_q (\sll ).$ On the other hand, the action of the subalgebra $U'_q (\sll )^{(2)}$  has level 0 and is completely reducible. The construction of the irreducible  decomposition of the level-1 $U_q (\sll )$-modules  relative to the level-0 action is the problem which we address in the present paper. To solve this problem we utilize, as our main tools, the semi-infinite wedges of \cite{KMS} and the Non-symmetric Macdonald polynomials  of \cite{Cheann}.   

As a result we obtain a parameterization of the irreducible components of the level-0 action by the skew Young diagrams of the border strip type, proving thereby a $q$-analogue of the conjecture made in \cite{KKN} in the classical setting. 

\section{The actions of the quantum affine algebra $U'_q(\sll )$}

\subsection{Definition of the quantum affine algebra $U'_q(\sll )$}

\begin{df} 
The quantum affine algebra $U'_q(\sll )$
is the unital
associative algebra over $\cplx $ with generators $E_i$, $F_i$, $K_i^{{}\pm
1}$\/ ($i\in I:=\{0,1, \dots ,n-1\}$) and the following defining relations:
\begin{eqnarray}K_iK_i^{- 1} =&1& =  K_i^{ -1}K_i,\\
K_iK_j&=&K_jK_i\;,\\
K_iE_jK_i^{-1}&=& q^{ a_{ij}}E_j \; ,\\
K_iF_jK_i^{-1}&=& q^{ -a_{ij}}F_j \; ,\\
{}[E_i , F_j]&=&\delta_{ij}\frac{K_i -K_i^{-1}}{q-q^{-1}}\;,
\end{eqnarray}
\vspace{-.15in}
\begin{equation}
\sum_{r=0}^{1-a_{ij}}
 (-1)^r\left[{1-a_{ij}\atop r}\right]_{q}\;
(E_i)^rE_j (E_i)^{1-a_{ij}-r}= 0\;,\;i\ne j .
\end{equation}
\vspace{-.2in}
\begin{equation}
\sum_{r=0}^{1-a_{ij}}
 (-1)^r\left[{1-a_{ij}\atop r}\right]_{q}\;
(F_i)^rF_j (F_i)^{1-a_{ij}-r}= 0\;,\;i\ne j .
\end{equation}
\vspace{-.15in}
\begin{equation}
 \mbox{where }
[n]_{q} :=\frac{q^n -q^{-n}}{q -q^{-1}},\; \;
{}\left[{n\atop r}\right]_{q}:=\frac{[n]_{q}[n-1]_{q}\ldots {}[n-r+1]_{q}}
{[r]_{q}[r-1]_{q}\ldots [1]_{q}},
\end{equation}
\vspace{-.2in}
\begin{equation}
 a_{ij}= \left\{
\begin{array}{cc}
2 & (i=j) \\
-1 &(|i-j|=1 ,(i,j)=(1,n) ,(n,1)) \\
0 & (\mbox{otherwise})
\end{array} \right.
 \; \; n \geq 3 ,
\end{equation}
\vspace{-.2in}
\begin{equation}
 a_{ij}= \left\{
\begin{array}{cc}
2 & (i=j) \\
-2 & (i\neq j) 
\end{array} \right.
 \; \; n =2 .
\end{equation}
\end{df}
The coproduct $\Delta $ is given by
\begin{eqnarray}
\Delta(E_{i})& =&  E_{i} \otimes K_{i} + 1 \otimes E_{i}, \label{cp1} \\
\Delta(F_{i})& =&  F_{i} \otimes 1 + K_{i}^{-1} \otimes F_{i}, \label{cp2} \\ 
\Delta(K_{i})& =&  K_{i} \otimes K_{i}.\label{cp3}
\end{eqnarray}

We put $c':=K_{0}K_{1} \dots K_{n-1}$ in $U'_{q}(\sll)$,
 then $c'$ is the central in $U'_{q}(\sll)$.

\subsection{q--wedge product and semi--infinite q--wedge product}

The affine Hecke algebra of type $\glN$, $\Hhat $ is a unital associative algebra over $\cplx[q^{\pm 1}]$ with generators $T_i^{\pm 1},$ $Y_j^{\pm 1},$ $i=1,2,\dots,N-1,$ $j=1,2,\dots
,N$ and relations 
\begin{gather*}
T_i T_i^{-1} =  T_i^{-1} T_i = 1, \quad (T_i + 1) (T_i - q^2) =0, \\
T_i T_{i+1} T_i =  T_{i+1} T_i T_{i+1},\\
T_i T_j = T_j T_i \quad \text{ if $ |j-i| > 1,$} \\
Y_i Y_j = Y_j Y_i, \quad T_i^{-1} Y_i T_i^{-1} = q^{-2} Y_{i+1}\\
Y_j T_i = T_i Y_j \quad \text{ if $ j \neq i,i+1,$}
\end{gather*}
The subalgebra $\HH$ generated by $ T_i^{\pm 1}$ is isomorphic to the Hecke algebra of type $\glN.$

Let $p \in \cplx^{\times}$ and consider the following operators in $\End(\cz)$ 
\begin{eqnarray}
& g_{i,j} = \frac{q^{-1}z_i - q z_j}{z_i - z_j}(K_{i,j} - 1) + q, \;\; & 1\leq i \neq j \leq N,\nonumber \\ 
& Y_i^{(N)} = g_{i,i+1}^{-1}K_{i,i+1} \cdots g_{i,N}^{-1} K_{i,N} p^{D_i} K_{1,i} g_{1,i} \cdots K_{i-1,i}g_{i-1,i}, \;\; & i=1,2,\dots,N, \nonumber
\end{eqnarray}
where $K_{i,j}$ acts on $\cz$ by permuting variables $z_i,z_j$ and $p^{D_i}$ is the difference operator  
\begin{equation}
p^{D_i}f(z_1,\dots,z_i,\dots,z_N)  =  f(z_1,\dots,pz_i,\dots,z_N) , \quad f \in \cplx[z_1^{\pm 1},\dots, z_N^{\pm 1}]. \nonumber 
\end{equation}
Then the assignment
\begin{equation}
T_i \mapsto \overset{c}{T_i} = -q g_{i,i+1}^{-1}, \quad Y_i \mapsto q^{1-N}Y_i^{(N)}
\end{equation}
defines a right action of $\Hhat $  on $\cz$. 

The commuting difference operators $Y_1^{(N)}, \dots, Y_N^{(N)}$ are called Cherednik's operators. 

Moreover, the assignment 
\begin{equation}
T_i \mapsto \overset{c}{T_i} = -q g_{i,i+1}^{-1}, \quad Y_i \mapsto z_i^{-1} \; \mbox{(multiplication)}
\end{equation}
defines another right action of $\Hhat $  on $\cz$. \\
\vspace{.1in}
\begin{rmk}
The actions of $-q g_{i,i+1}^{-1}, \; q^{1-N}Y_i^{(N)}, \; z_i^{-1}$ are related to the toroidal Hecke algebra \cite{VV} or the double affine Hecke algebra \cite{Cheann}.
\end{rmk}
\vspace{.1in}

Let $V = \cplxn$, with basis $\{ v_1,\dots,v_n \}.$ Then $\otimes^N V$ admits a left $\HH$-action given by    
\begin{gather}
T_i \mapsto \overset{s}{T_i} = 1^{\otimes^{i-1}}\otimes \overset{s}{T} \otimes 1^{\otimes^{N-i-1}}, \quad \text{ where } \quad  \overset{s}{T}  \in \End(\otimes^2 V) \\
\text{ and } \quad   \overset{s}{T}(v_{\ep_1}\otimes v_{\ep_2}) = \begin{cases} q^2 v_{\ep_1}\otimes v_{\ep_2} & \text{ if $ \ep_1 = \ep_2,$} \\
 q v_{\ep_2}\otimes v_{\ep_1} & \text{ if $ \ep_1 < \ep_2,$ }\\
q v_{\ep_2}\otimes v_{\ep_1} + (q^2  - 1) v_{\ep_1}\otimes v_{\ep_2}& \text{ if $ \ep_1 > \ep_2.$ }  \end{cases}
\end{gather}

Let $V(z)= \cplx [z^{\pm 1}]\otimes V,$ with  basis $ \{ z^m \otimes v_{\ep} \}$, $ m\in \zint $ , $ \ep \in \{1,2,\dots,n\} $. Often it will be convenient to set $ k = \ep - nm $ and $ u_k = z^m \otimes v_{\ep} $.
Then $\{ u_k \}, $  $k \in \zint $ is a basis of $V(z)$. In what follows we will write $ z^m v_{\ep} $ as a short-hand for $ z^m \otimes v_{\ep} ,$ and use both  notations: $ u_k $ and $ z^m  v_{\ep} $ switching between them according to convenience. 
The two actions of the Hecke algebra are naturally extended on the tensor product $\cz\otimes(\otimes^N V)$ so that $\overset{c}T_i$ acts trivially on $\otimes^N V$ and  $\overset{s}T_i$ acts trivially on $\cz.$ The vector space $\otimes^N V(z)$ is identified with  $\cz\otimes(\otimes^N V)$ and the $q$-wedge product \cite{KMS} is defined as the following quotient space:
\begin{equation}
\wedge^N V(z) = \otimes^N V(z)/ \sum_{i=1}^{N-1}\Ker\left( \overset{c}T_i + q^2 (\overset{s}T_i)^{-1} \right). \label{e: q-wedgeprod}
\end{equation}

Let $\Lambda: \otimes^N V(z) \rightarrow \wedge^N V(z) $ be the quotient map specified by (\ref{e: q-wedgeprod}). The image of a pure tensor $ u_{k_1}\otimes u_{k_2}\otimes \dots \otimes u_{k_N} $ under this map is called a wedge and is denoted by
\begin{equation} 
u_{k_1}\wedge u_{k_2}\wedge \dots \wedge u_{k_N} :=\Lambda ( u_{k_1}\otimes u_{k_2}\otimes \dots \otimes u_{k_N}) . \label{e: wedge} 
\end{equation}
A wedge is normally ordered if $k_1 > k_2 > \cdots > k_N.$ In \cite{KMS} it is proven that normally ordered wedges form a basis in $\wedge^N V(z)$. 
\vspace{.2in}

Let us now define the semi-infinite q--wedge product $\wedge ^{\frac{\infty}{2}} V(z)$ and for any integer $M$ its subspace $F_M$, following \cite{KMS}.

Let $\otimes ^{\frac{\infty}{2}} V(z)$ be the space spanned by the vectors
$u_{k_{1} } \otimes u_{k_{2}} \otimes \dots , \; ( k_{i+1} = k_{i} -1 , \; i >>1 )$. We define the space $\wedge ^{\frac{\infty}{2}} V(z)$ as the quotient of $\otimes ^{\frac{\infty}{2}} V(z)$:
\begin{equation}
\wedge ^{\frac{\infty}{2}} V(z) :=  \otimes ^{\frac{\infty}{2}} V(z) / \sum_{i=1}^{\infty }\Ker\left( \overset{c}T_i + q^2 (\overset{s}T_i)^{-1}\right). \label{e: inf-wedgeprod}
\end{equation}
Let $\Lambda: \otimes ^{\frac{\infty}{2}} V(z) \rightarrow \wedge ^{\frac{\infty}{2}} V(z)$ be the quotient map specified by (\ref{e: inf-wedgeprod}). The image of a pure tensor $ u_{k_1}\otimes u_{k_2}\otimes \dots  $ under this map is called a semi-infinite wedge and is denoted by
\begin{equation} 
u_{k_1}\wedge u_{k_2}\wedge \dots  :=\Lambda ( u_{k_1}\otimes u_{k_2}\otimes \dots ) . \label{e: infwedge} 
\end{equation}
A semi-infinite wedge is normally ordered if $k_1 > k_2 > \cdots $ and $k_{i+1}= k_i -1 \; (i >> 1)$. In \cite{KMS} it is proven that normally ordered semi-infinite wedges form a basis in $\wedge ^{\frac{\infty}{2}} V(z)$.
\vspace{.1in}

Let $U_{M}$ be the subspace of $\otimes ^{\frac{\infty}{2}} V(z)$ spanned by the vectors $u_{k_{1} } \otimes u_{k_{2}} \otimes \dots , \; ( k_{i} = M-i+1 , \; i >>1 )$.
Let $F_{M}$ be the quotient space of $U_{M}$ defined by the map (\ref{e: infwedge}).
 Then $F_{M}$ is a subspace of $\wedge ^{\frac{\infty}{2}} V(z)$, and the vectors $u_{k_1}\wedge u_{k_2}\wedge \dots $, ($k_1 > k_2 > \cdots $,  $k_{i} = M-i+1, \; i >>1 )$ form a basis of $F_{M}$. We will  call the space $F_{M}$ the q--deformed Fock space.

\subsection{Actions of the quantum affine algebra on the q--wedge product}
\label{fwrep}

We will define two actions of $U'_q (\sll )$ on the space $\wedge ^N V(z)$.

The first one is defined as follows.
\begin{align}
& E_i ( m \otimes v)  = \sum_{j=1}^N m \otimes E_j^{i,i+1} K_{j+1}^{i}\dots K_N^{i} v,  \label{e: Efin}\\  
& F_i ( m \otimes v)  = \sum_{j=1}^N m \otimes (K_{1}^{i})^{-1}\dots (K_{j-1}^{i})^{-1}E_j^{i+1,i} v,  \label{e: Ffin}\\  
& K_i( m \otimes v)   = m \otimes   K^{i}_1K^{i}_2\dots K^{i}_N v ,\quad    (i=1,2,\dots,n-1)  \label{e: Kfin} \\ 
& E_0( m \otimes v) =  \sum_{j=1}^N m Y_j^{-1}\otimes E_j^{n,1} K_{j+1}^{0}\dots K_N^{0} v,  \label{e: E0} \\ 
& F_0 ( m \otimes v)    =   \sum_{j=1}^N m Y_j \otimes (K_{1}^{0})^{-1}\dots (K_{j-1}^{0})^{-1}E_j^{1,n} v,  \label{e: F0} \\  
&  K_0  = (K_1 K_2 \cdots K_{n-1})^{-1}. \label{e: K0}
\end{align}
Here $ E_j^{i,k} = 1^{\otimes^{j-1}}\otimes E^{i,k} \otimes 1^{\otimes^{N-j}},$ where $E^{i,k} \in \End(V)$ is the matrix unit in the basis $ v_1,\dots,v_n,$ and $ K_j^i = q^{E_j^{i,i} - E_j^{i+1,i+1}},$ $ K_j^0 = (K_j^1 K_j^2 \cdots K_j^{n-1})^{-1}.$ 

We will denote this action by $U_0^{(N)}$. Note that  it is well defined on the quotient space $\wedge ^N V(z)$ in view of the relations of the affine Hecke algebra.

The second one is defined as follows.
\begin{align}
 E_0( m \otimes v) =  \sum_{j=1}^N m z_j \otimes E_j^{n,1} K_{j+1}^{0}\dots K_N^{0} v,   \\ 
 F_0 ( m \otimes v)    =   \sum_{j=1}^N m z_j^{-1} \otimes (K_{1}^{0})^{-1}\dots (K_{j-1}^{0})^{-1}E_j^{1,n} v. 
\end{align}
The actions of other generators are the same as in (\ref{e: Efin}--\ref{e: Kfin}, \ref{e: K0}).

We will denote this action by $U_1^{(N)}$. Again, this action is well defined on the quotient space $\wedge ^N V(z)$ in view of the relations of the affine Hecke algebra.

\subsection{Level--0 action of the quantum affine algebra on the q--deformed Fock space} \label{sec:L-0}

We will define a level-0 action of $U_q'(\sll)$ on $F_M$ $(M\in \zint)$ following the paper \cite{TU,STU}.

Let $ \ee := (\ep_1,\ep_2,\dots,\ep_N) $ where $ \ep_i \in \{1,2,\dots,n\} $. For a sequence $\ee $ we set
\begin{equation}
{\bold v}_{\ee}:= v_{\ep_1}\otimes v_{\ep_2}\otimes \dots \otimes v_{\ep_N}\quad  ( \in \otimes^N \cplxn ).  
\end{equation}
A sequence $\mm := (m_1,m_2,\dots,m_N)$ from $\zint^N$ is called $n$-strict if it contains no more than $n$ equal elements of any given value. 
Let us define the sets $\MC_N^n$ and $\EC(\mm)$ by 
\begin{align}
 & \MC_N^n := \{\mm = (m_1,m_2,\dots,m_N) \in \zint^N \; | \; m_1\leq m_2 \leq \dots \leq m_N,\; \text{$\mm $  is $n$-strict } \}, \label{e: Mm}\\
 \intertext{ and for $ \mm \in \MC_N^n $ } 
& \EC(\mm) := \{ \ee = (\ep_1,\ep_2,\dots,\ep_N) \in \{1,2,\dots,n\}^N \; | \; \ep_i > \ep_{i+1} \; \text{for all $i$ s.t. $m_i =m_{i+1}$ } \}.   \label{e: Em}
\end{align}
In these notations the set 
\begin{equation}
\{ w(\mm,\ee):= \Lambda( \zz^{\mm} \otimes {\bold v}_{\ee}) = z^{m_1}v_{\ep_1}\wedge z^{m_2}v_{\ep_2}\wedge\dots\wedge z^{m_N}v_{\ep_N} \quad | \quad \mm \in \MC_N^n , \ee \in \EC(\mm) \}.
\end{equation}
is nothing but the base of the normally ordered wedges in $\wedge^N V(z)$. We  will use the notation $w(\mm,\ee)$ {\em exclusively} for normally ordered wedges.

Similarly for a semi-infinite wedge $w = u_{k_1}\wedge u_{k_2} \wedge \dots \; $ $=$ $z^{m_1}v_{\ep_1}\wedge z^{m_2}v_{\ep_2} \wedge \dots \;,$ such that $w \in F_M,$ the semi-infinite sequences $\mm = (m_1,m_2,\dots \;\;)$ and $\ee = (\ep_1,\ep_2,\dots \;\;)$ are defined by $ k_i = \ep_i - nm_i $, $\ep_i \in \{1,2,\dots,n\},$ $m_i \in \zint.$ In particular the $\mm$- and $\ee$- sequences of the vacuum vector in $F_M$ will be denoted by $\mm^0$ and $\ee^0$:      
\begin{equation}
 |M\rangle = u_M\wedge u_{M-1} \wedge u_{M-2}\wedge  \dots \; = z^{m_1^0}v_{\ep_1^0}\wedge z^{m_2^0}v_{\ep_2^0} \wedge z^{m_3^0}v_{\ep_3^0} \wedge \cdots \;. \end{equation}

The Fock space $F_M$ is $\zintp$-graded. For any semi-infinite wedge $w$ $=$ $u_{k_1}\wedge u_{k_2} \wedge \dots \;$ $=$  $z^{m_1}v_{\ep_1}\wedge z^{m_2}v_{\ep_2} \wedge \dots \;$ $\in F_M$ the degree $| w |$ is defined by  
\begin{equation}
| w | = \sum_{i \geq 1} m_i^0 - m_i . \label{wdeg}
\end{equation}
Let us denote by $F_M^k$ $\subset$ $F_M$ the homogeneous component of degree $k$.

We will define a level--0 action of $U_q'(\sll)$ on the Fock space $F_M$ in such a way that each homogeneous component $F_M^k$ will be invariant with respect to this action. Throughout this section we fix an integer $M$ and  $s \in \{0,1,2,\dots,n-1 \} $ such that  $M = s\bmod n.$ 

Let $l$ be a non-negative integer and define $ V_M^{s+nl} \subset \wedge^{s+nl}V(z) $ as follows:
\begin{equation}
V_M^{s+nl}  = \bigoplus\begin{Sb} \mm \in \MC_{s+nl}^n, \ee \in \EC(\mm) \\ m_{s+nl} \leq m_{s+nl}^0 \end{Sb} \cplx w(\mm,\ee) . \label{e: VM}
\end{equation}

The vector space $V_M^{s+nl}$ has a grading similar to the grading of the Fock space $F_M$. In this case the degree $|w |$ of a wedge  $w$ $=$ $u_{k_1}\wedge u_{k_2} \wedge \dots \wedge u_{k_{s+nl}}$ $=$  $z^{m_1}v_{\ep_1}\wedge z^{m_2}v_{\ep_2} \wedge \dots \wedge z^{m_{s+nl}}v_{\ep_{s+nl}} $ $\in$  $V_M^{s+nl}$ is defined by     
\begin{equation}
| w | = \sum_{i = 1}^{s+nl} m_i^0 - m_i .
\label{degifw}
\end{equation}
The degree is a non-negative integer, and for $k\in \zintp$ we denote by $ V_M^{s+nl,k} $ the homogeneous component of degree $k.$

The following result is contained in the paper \cite{TU}:
\begin{prop}
For each $k\in \zintp$ the homogeneous component $ V_M^{s+nl,k} $ $\subset$ $\wedge^{s+nl}V(z)$ is invariant under the $U_q'(\sll) $-action  $U_0^{(s+nl)}$ defined in section $\ref{fwrep}$. 
\end{prop}
We have $|\oro(w)|$ $=$ $|w|$ and hence $\oro$ $:$ $V_M^{s+nl,k} $ $\rightarrow$ $F_M^k$ for all $k \in \zintp$. In the paper \cite{STU} the following propositions are shown.
\begin{prop}
When $l \geq k$ the map $\oro$ is an isomorphism of vector spaces. \label{p: isom}
\end{prop}
\begin{prop} \label{u1isom}
For each triple of non-negative integers $k,l,m$ such that $k\leq l < m$ the map $\orom$ $:$ $ V_M^{s+nl,k}$ $\rightarrow$ $ V_M^{s+nm,k},$ defined for any $w$ $\in$  $ V_M^{s+nl,k}$ by 
\begin{equation}
\orom(w) = w\wedge u_{M-s-nl}\wedge u_{M-s-nl-1}\wedge \dots \wedge u_{M-s-nm+1},
\end{equation}
is an isomorphism of the $U_q'(\sll)$-modules. 
\end{prop}
We define on the vector space $F_M^k$ a level--0 action of $U_q'(\sll)$ by using Propositions \ref{p: isom} and \ref{u1isom}:
\begin{df} 
The vector space $F_M^k$ is a level-$0$ module of $U_q'(\sll)$ with the action $U_0$ defined by  
\begin{equation}
 U_0 = \oro U_0^{(s+nl)} {\oro}^{-1} \quad \text{ where $ l \geq k .$ }  \label{eq:U0def}
\end{equation} 
\end{df}
This definition does not depend on the choice of $l$ as long as $l$ is greater or equal to $k.$ Since we have
\begin{equation}
F_M = \bigoplus_{k \geq 0} F_{M}^k
\end{equation}
the level-0 action $U_0$ extends to the entire Fock space $F_M.$

\subsection{Level--1 action of the quantum affine algebra on the q--deformed Fock space}

In this section we review the level--1 action of $U_q'(\sll)$ on the Fock space $F_M$ \cite{KMS}. 

First we define the action of $U'_q(\sll)$ (generated by $E _{i}, \; F _i , \; K_i, \; $ $i=0, \dots ,n-1$) on the vector $|M' \rangle $ as follows.
\begin{equation}
E_i| M' \rangle =0,
\label{eq79}
\end{equation}
\vspace{-.18in}
\begin{equation}
F_i| M' \rangle = \left\{
\begin{array}{ll}
u_{M' +1}\wedge u_{M'-1}\wedge u_{M'-2}\wedge \cdots &  \mbox { if } i\equiv M' \mbox{ mod }n; \\
0 & \mbox{ otherwise},
\end{array}
\right.
\end{equation}
\vspace{-.18in}
\begin{equation}
K_i | M' \rangle = \left\{
\begin{array}{ll}
q| M' \rangle & \mbox { if } i\equiv M' \mbox{ mod }n; \\
| M' \rangle & \mbox{ otherwise},
\end{array}
\right.
\end{equation}
For every element $v \in F_{M}$, there exists $N$ such that
\begin{equation}
v = v^{(N)} \wedge | M -N \rangle , \; \; \; \; \; 
v^{(N)} \in \wedge ^N V(z).
\end{equation}
We define the actions of $E _{i}, \; F _i , \; K_i, \; $ $i=0, \dots ,n-1$ on the vector $v$ as follows.
\begin{eqnarray}
& E _i v := E _i v^{(N)} \wedge K_i | M -N \rangle +  v^{(N)} \wedge E_i | M -N
\rangle ,& \\
& F _i v := F _i v^{(N)} \wedge | M -N \rangle +  K_i^{-1} v^{(N)} \wedge F_i  
M -N \rangle ,& \\
& K _i v := K _i  v^{(N)} \wedge K _i | M -N \rangle. &
\end{eqnarray}
The actions of $E _{i}, \; F _i , \; K_i, \; $ $i=0, \dots ,n-1$ on $v^{(N)}$ are determined in Section \ref{fwrep}. The definition of the actions on 
$v$ does not depend on $N$ and is well--defined, and we can easily check that the $U'_q (\sll)$-module defined in this section is level-1. We will use the notation $U_1$ for this  $U'_q (\sll)$-action on the Fock space. \\
\begin{rmk}
The two actions $U_0$ and $U_1$ appear as the representations of the subalgebras of the quantum toroidal algebra. For details, see \cite{STU}.
\end{rmk}

\subsection{The $p=1$ case}

In the paper \cite{KMS} it was demonstrated that the Fock space $F_M$ admits an action of the Heisenberg algebra $H$ which commutes with the level-1 action  $U_1$ of the algebra $U_q'(\sll).$ The Heisenberg algebra is a unital $\cplx$-algebra generated by elements $1, B_a$ with $a \in \zint_{\neq 0}$ which are subject to  relations 
\begin{equation}
[ B_a , B_b ] = \delta_{a+b,0}a\frac{1 - q^{2n a}}{1 - q^{2 a}} . \label{e: Heisenbergrel} \end{equation}
The Fock space $F_M$ is an $H$-module with the action of the generators given by \cite{KMS} 
\begin{equation}
B_a = \sum_{i=1}^{\infty} z_i^{a}. \label{e: Bact}
\end{equation}

Let $\cplx[H_{-}]$ be the Fock space of $H$, i.e.,  $\cplx[H_{-}]= \cplx[B_{-1},B_{-2},\dots,\:].$ The element $B_{-a}$ $ ( a = 1,2,\dots\:)$ acts on  $\cplx[H_{-}]$ by multiplication. The action of  $B_{a}$ $ ( a = 1,2,\dots\:)$ is given by (\ref{e: Heisenbergrel}) together with the relation
\begin{equation}
 B_a \cdot 1 = 0 \qquad \text{ for $ a \geq 1. $} 
\end{equation}

Let $\Lambda_i$ $(i\in \{ 0,1,\dots,n-1\}) $ be the fundamental weights of $\sll'.$ And let $V(\Lambda_i)$ be the irreducible (level-1) highest weight module of $U_q'(\sll)$ with highest weight vector $V_{\Lambda_i}$ and highest weight $\Lambda_i.$

The following results are proven in \cite{KMS}:
\begin{itemize}
\item The action of the Heisenberg algebra on $F_M$ and the action $U_1$ of $U'_q (\sll )$ commute. 
\item There is an isomorphism 
\begin{equation}
 \iota_M : F_M \cong V(\Lambda_i)\otimes \cplx[H_{-}] \qquad ( M= i\bmod n )
\label{e: iota} \end{equation}
of $U_q'(\sll)\otimes H$-modules normalized so that $\iota_M( |M\rangle )$ $ =$ $V(\Lambda_i)\otimes 1.$
\end{itemize}
In general the level-0 $U_q'(\sll)$-action $U_0$   does not commute with the Heisenberg algebra. However if we choose the parameter $p$ in $U_0$ in a special way, then $U_0$ commute with the negative frequency part of $H$. Precisely, we have the following proposition, proved in \cite{STU}:  
\begin{prop} \label{p: U0comH-}
At $p=1$ we have
\begin{equation}
[ U_0 , H_{-}] = 0.
\end{equation}
\end{prop}

Let $H_{-}'$ be the non-unital subalgebra in $H$ generated by  $ B_{-1},B_{-2},\dots\;.$ Proposition \ref{p: U0comH-} allows us to define a level-0 $U_q'(\sll)$-module structure on the irreducible level-1 module $V(\Lambda_i)$ $(i \in \{0,1,\dots,n-1\}).$ Indeed from this proposition it follows that the subspace  
\begin{equation}
  H_{-}'F_M \subset F_M 
\end{equation}
is invariant with respect to the action $U_0$ at $p=1$ and therefore a level-0 action of $U_q'(\sll)$ is defined on the quotient space
\begin{equation}
F_M/ (H_{-}'F_M)  \label{eq: jjj}
\end{equation}
which in view of (\ref{e: iota}) is isomorphic to  $V(\Lambda_i)$ with $ i \equiv M$  mod $n$.

\section{Skew Young diagrams and the level--$0$ representations of $U'_q(\sll )$}

\subsection{Skew Young diagrams}
Let us  recall, following the book \cite{Macbook}, the definitions of the skew (Young) diagrams, their semi-standard   tableaux and the associated skew Schur functions.

Let $\lambda ,\mu $ be partitions i.e. sequences of non--negative integers.
We assume $\l _{i} \geq \mu _i $ for all possible $i$, and if $\mu _j < i \leq \l_j $ then we draw a square whose edges are $(i-1,j-1), \; (i-1, j), \; (i,j)$ and $(i,j-1)$. (For example, see Figure 1.) This diagram is called a skew (Young) diagram and is denoted as  $\l \setminus \mu $.
We define the degree of the skew Young diagram $\l \setminus \mu $ as $| \l \setminus \mu |$ = $\sum _{i} (\l _{i} - \mu _i)$.

A skew diagram is called a border strip if it is connected and contains no $2 \times 2$ blocks of boxes. Let $\langle m_1 , \dots , m_r \rangle $ denote the border strip of $r$ columns such that the length of $i$--th column (from the right) is $m_i$. (Figure 1)

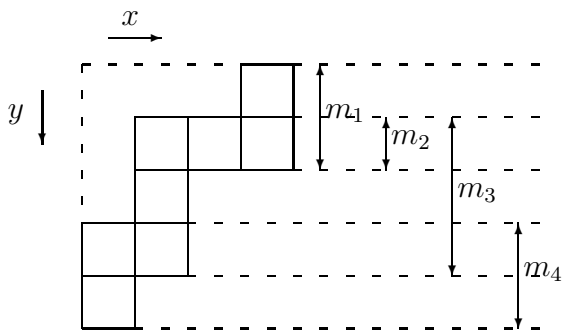
\begin{figure} 
\begin{center}
\begin{picture}(180,140)(0,-130)
\put(60,0){\line(1,0){20}}
\put(20,-20){\line(1,0){60}}
\put(20,-40){\line(1,0){60}}
\put(0,-60){\line(1,0){40}}
\put(0,-80){\line(1,0){40}}
\put(0,-100){\line(1,0){20}}
\put(0,-60){\line(0,-1){40}}
\put(20,-20){\line(0,-1){80}}
\put(40,-20){\line(0,-1){60}}
\put(60,0){\line(0,-1){40}}
\put(80,0){\line(0,-1){40}}
\multiput(0,0)(10,0){6}{\line(1,0){3}}
\multiput(0,0)(0,-10){6}{\line(0,-1){3}}
\multiput(80,0)(10,0){10}{\line(1,0){3}}
\multiput(80,-20)(10,0){10}{\line(1,0){3}}
\multiput(80,-40)(10,0){10}{\line(1,0){3}}
\multiput(40,-60)(10,0){14}{\line(1,0){3}}
\multiput(40,-80)(10,0){14}{\line(1,0){3}}
\multiput(20,-100)(10,0){16}{\line(1,0){3}}
\put(90,-20){\vector(0,1){20}}
\put(90,-20){\vector(0,-1){20}}
\put(92,-20){$m_1$}
\put(115,-30){\vector(0,1){10}}
\put(115,-30){\vector(0,-1){10}}
\put(117,-30){$m_2$}
\put(140,-50){\vector(0,1){30}}
\put(140,-50){\vector(0,-1){30}}
\put(142,-50){$m_3$}
\put(165,-80){\vector(0,1){20}}
\put(165,-80){\vector(0,-1){20}}
\put(167,-80){$m_4$}
\put(-15,-10){\vector(0,-1){20}}
\put(-28,-20){$y$}
\put(10,10){\vector(1,0){20}}
\put(15,15){$x$}
\end{picture}
\end{center}
\caption{$\lambda = (4,4,2,2,1), \; \mu= (3,1,1) $.}
\end{figure}

A semi--standard tableau (s.s.t.) of the skew diagram $\l \setminus \mu $ is obtained by inscribing  integers $1,2 ,\dots ,n$ in  each square of the skew diagram. The rule of the semi--standard tableau is as follows. The numbers are strictly increasing along the column and weakly increasing along the row.
For each semi--standard tableau $T$, let  $n_i(T)$ be the multiplicity  of $i$ in $T.$ 
\begin{df}
For each skew diagram $ \l \setminus \mu $, the skew Schur function $s_{\l \setminus \mu}$ is defined as follows. 
\begin{equation}
s_{\l \setminus \mu} (z)= \sum_{T} z_1^{n_{1}(T)} z_2^{n_{2}(T)} \dots z_N^{n_{N}(T)} .
\end{equation}
Here the summmation  is over the set of semi--standard tableaux of the skew diagram $\l \setminus \mu $. 
\end{df}

\subsection{The level--$0$ representations of $U'_q(\sll )$ associated with the skew diagrams} 

Fix a skew diagram $\l \setminus \mu $ of the border strip type and degree $N$. We put a number $(1, \dots ,N)$ on each box such that if $l>k$ then $x_l >x_k$ or ($x_l = x_k$ and $y_l > y_k$), where $(x,y)$ is a box contained in the skew diagram, and set $a_l = -2x_l + 2y_l+a$ ($a$ is fixed). (Figure 2)

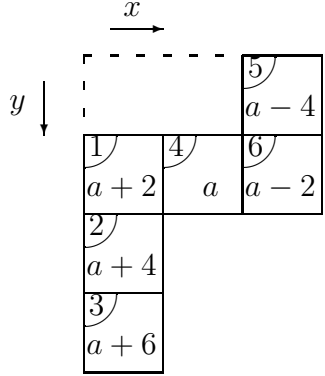
\begin{figure} 
\begin{center}
\begin{picture}(100,140)(0,-130)
\put(60,0){\line(1,0){30}}
\put(0,-30){\line(1,0){90}}
\put(0,-60){\line(1,0){90}}
\put(0,-90){\line(1,0){30}}
\put(0,-120){\line(1,0){30}}
\put(0,-30){\line(0,-1){90}}
\put(30,-30){\line(0,-1){90}}
\put(60,0){\line(0,-1){60}}
\put(90,0){\line(0,-1){60}}
\multiput(0,0)(10,0){6}{\line(1,0){3}}
\multiput(0,0)(0,-10){3}{\line(0,-1){3}}
\put(1,-53){$a+2$}
\put(1,-83){$a+4$}
\put(1,-113){$a+6$}
\put(45,-53){$a$}
\put(61,-23){$a-4$}
\put(61,-53){$a-2$}
\put(0,-30){\oval(25,25)[br]}
\put(0,-60){\oval(25,25)[br]}
\put(0,-90){\oval(25,25)[br]}
\put(30,-30){\oval(25,25)[br]}
\put(60,0){\oval(25,25)[br]}
\put(60,-30){\oval(25,25)[br]}
\put(2,-38){$1$}
\put(2,-68){$2$}
\put(2,-98){$3$}
\put(32,-38){$4$}
\put(62,-8){$5$}
\put(62,-38){$6$}
\put(-15,-10){\vector(0,-1){20}}
\put(-28,-20){$y$}
\put(10,10){\vector(1,0){20}}
\put(15,15){$x$}

\end{picture}
\end{center}
\caption{The case $\langle 2,1,3 \rangle $.}
\end{figure}

On the space $\otimes ^{N} V$, we define the  evaluation action $ \pi_{a_1,\dots,a_N}^{(N)} $ of $U_q'(\sll)$ 
\begin{align}
& \pi_{a_1,\dots,a_N}^{(N)}( E_i ) = \sum_{j=1}^N E_j^{i,i+1} K_{j+1}^{i}\dots K_N^{i} ,\\  
& \pi_{a_1,\dots,a_N}^{(N)}( F_i ) = \sum_{j=1}^N (K_{1}^{i})^{-1}\dots (K_{j-1}^{i})^{-1}E_j^{i+1,i} ,\\  
& \pi_{a_1,\dots,a_N}^{(N)}( K_i ) =  K^{i}_1K^{i}_2\dots K^{i}_N , \quad (i=1,2,\dots,n-1)  \\ 
& \pi_{a_1,\dots,a_N}^{(N)}( E_0 ) = \sum_{j=1}^N q^{a_j} E_j^{n,1} K_{j+1}^{0}\dots K_N^{0} ,\\ 
& \pi_{a_1,\dots,a_N}^{(N)}( F_0 ) = \sum_{j=1}^N q^{-a_j} (K_{1}^{0})^{-1}\dots (K_{j-1}^{0})^{-1}E_j^{1,n} , \\  
& \pi_{a_1,\dots,a_N}^{(N)}( K_0 ) = (\pi_{a_1,\dots,a_N}^{(N)}(K_1 K_2 \cdots K_{n-1}))^{-1},
\end{align}
and on the same space, we consider the following operators:
\begin{equation}
R_{i,j} (x) = \frac{x S_{i,j}^{-1}- S_{i,j}}{x-1} P_{i,j}, \; \; 
\check{R}_{i,j} (x) = \frac{x S_{i,j}^{-1}- S_{i,j}}{x-1},
\label{Rcheck}
\end{equation}
where $P_{i,j}(\dots \otimes \stackrel{i}{u} \otimes \dots \otimes \stackrel{j}{v} \otimes \dots ) = \dots \otimes \stackrel{i}{v} \otimes \dots \otimes \stackrel{j}{u} \otimes \dots$.
We define 
\begin{align}
& R_{\l \setminus \mu } = \prod _{1\leq i< j \leq N} R_{i,j} (q^{a_i-a_j}), \\
& \bar{R}_{\l \setminus \mu } = \prod _{1\leq i< j \leq N} R_{j,i} (q^{a_i-a_j}), \\
& \check{R}_{\l \setminus \mu } = \prod _{1\leq i< j \leq N} \check{R}_{N+i-j,N+i-j+1} (q^{a_i-a_j}), 
\end{align}
where $(i,j) $ is on the right to $(i',j')$ in the product if $i<i'$ or ($j<j'$ and $i=i'$). As a special case of \cite{CheDuke} Proposition 1.5., we have 
\begin{prop}[\cite{CheDuke}] \label{skewrep}
The subspace $\mbox{Im} R_{\l \setminus \mu }(\otimes ^N V) (= \mbox{Im} \check{R}_{\l \setminus \mu }(\otimes ^N V)) $ with the action $\pi_{a_1,\dots,a_N}^{(N)} $ is an irreducible $U_q '(\sll )$--module, and the map $\check{R}_{\l \setminus \mu }$ : $(\pi_{a_N,\dots,a_1}^{(N)} , \; \otimes ^N V / $ Ker $\bar{R}_{\l \setminus \mu })=  (\pi_{a_N,\dots,a_1}^{(N)} , \;\otimes ^N V / $ Ker $\check{R}_{\l \setminus \mu }))$ $\rightarrow $ $(\pi_{a_1,\dots,a_N}^{(N)} , \; $ Im $\check{R}_{\l \setminus \mu })$ is an isomorphism of the $U_q '(\sll )$--modules.
\end{prop} 
\begin{rmk}
In \cite{CheDuke}, this proposition is proved in the $U_q '(\widehat{ {\frak{g}} {\frak{l}}}_n)$--module case and the normalizations of $q$ and $x$ are different. The irreducibility as the $U_q '(\sll )$--module follows from the result of \cite{AK}.
\end{rmk}

\subsection{Character formulas}
Let $\overline{\Lambda }_i $ $(i=1, \dots ,n-1)$ be the fundamental weights of $\sln $ and let $\epsilon _i = \overline{\Lambda }_i - \overline{\Lambda }_{i-1} $ $(i=1, \dots ,n)$ with $\overline{\Lambda }_0 = \overline{\Lambda }_n :=0$.

The subalgebra  of $U'_q (\sll )$ generated by $E_i , \; F_i ,\; K_i^{\pm }$ $(i=1,\dots n-1)$ is isomorphic to the algebra $U_q (\sln)$.  In the paper \cite{KKN} the $\sln$--character of the irreducible $Y(\sln )$--representation associated with a skew diagram was shown to be given by the corresponding skew Schur function. This result  is immediately generalized to the $q$-deformed situation. Precisely we have 

\begin{prop}[\cite{KKN}] The skew Schur function $s_{\l \setminus \mu }(z)$ where $z_i = e ^{\epsilon _i}$ is equal to the 
 $U_q (\sln)$-character of the irreducible $U_q '(\sll )$--module described by Proposition \ref{skewrep}.
\end{prop} 
As a corollary we get
\begin{cor} \label{dimbs}
The dimension of the space $\mbox{Im} R_{\l \setminus \mu } \subset (\otimes ^N V)$ is equal to the total number of the semi--standard tableaux of the skew diagram $\l \setminus \mu $.
\end{cor}

Let $V(\Lambda _k )$ be the level--$1$ irreducible module of $U_q (\sll)$ whose highest weight is the $k$--th fundamental weight $\Lambda _k$ of $\sll '$. We set $ch (V(\Lambda _k )) = \sum _{i,\l } ($dim$V_{\l , i}) e ^{\l} q^{i}$, where $V_{\l ,i}$ is the weight subspace with  $U_q (\sln)$-weight $\l$ and homogeneous degree $i.$
The following proposition is proved in \cite{KKN}.
\begin{prop} [\cite{KKN}] \label{chara}
Setting  $z_i = e ^{\epsilon _i}$ we have
\begin{equation}
ch (V(\Lambda _k )) =   q^{\frac{1-n}{24}- \frac{k(n-k)}{2n}} \sum_{\theta \in BS \atop {| \theta | \equiv k \; mod \; n}} q ^{\frac{1}{2n} | \theta | (n- | \theta |)+t(\theta )} s_\theta (z).
\end{equation}
where BS is the set of all the border strips $\theta = \langle m_1 , \dots , m_r \rangle $ and $t(\theta ) = \sum_{i=1}^{r-1} (r-i)m_i$ with $m_r <n$.
\end{prop}
Note that if $m_i >n$ for some $i$, then the skew Schur function $s_{\theta }$ is equal to $0$, moreover, for  the border strip of the form $\theta _l = \langle m_1, \dots , m_r , \overbrace{n, \dots n}^{l} \rangle $ the number $\frac{1}{2n} | \theta _l | (n- | \theta _l |)+t(\theta_l  )$ does not depend on $l$.

\section{Non--symmetric Macdonald polynomials and the decomposition}

\subsection{Non--symmetric Macdonald polynomials}
\label{nsMac}

We will define the non--symmetric Macdonald polynomials as the joint eigenfunctions of the Cherednik's operators $Y_i^{(N)}$ $(i=1 , \dots , N)$ \cite{U1,TU}. It will be convenient for our  purposes to label these polynomials by the set of pairs $(\l , \s)$ which we now describe.

Let $\MN $ be the a set of all non--decreasing sequences of integers  
$ \l =  (\lambda _{1} \: , \: \lambda _{2} \: , \dots , \:  \lambda _{N}), $ 
and  let  $\MN ^n $ be the subset of $\MN $ which consists of all  $n$--strict non--decreasing sequences (cf. Section  \ref{sec:L-0}). 
For each $\lambda \in \MN $ we set $|\lambda |$ :=$ \sum_{i=1}^{N} \lambda _{i}. $ 
 For $\lambda , \mu \in \MN $ such that $|\lambda|=|\mu|$  we define the dominance (partial) ordering:
\begin{equation}
 \lambda \succeq \mu \hspace{.4in} \Leftrightarrow 
\hspace{.4in} \sum_{j=1}^{i} \lambda_{j} \geq \sum_{j=1}^{i} \mu_{j} \; \;
 (i\in \setN). 
\end{equation} 
 Let $S^{\lambda } \subset {{\frak{S}}_{N}} $ be the set of
 elements ${\sigma }$ such that
 if $\lambda_{\sigma (i)} = \lambda_{\sigma (j)}$ and $\sigma (i) < \sigma (j) $ then $i<j$.
 We define the total ordering on  $S^{\lambda }$:   
\begin{equation}
 \sigma \succ \sigma ' \; \; \Leftrightarrow \; \; 
\mbox{the last nonzero element of } 
 (\lambda _{\sigma (i)} -\lambda _{\sigma' (i) } )_{i=1}^N   \mbox{ is } < 0. 
\label{Sord}
\end{equation}
Then the following properties are satisfied. (In what follows the  $\s (i,i+1) $ denotes  the composition of $\sigma $ and a transposition $(i,i+1)$)\\
a) $S^{\lambda }$ has the unique minimal element with respect to the ordering (\ref{Sord}). We denote this element  by $min.$ Note that we one has $\l _{min(i)} \leq \l _{min(i+1)} $ $(i=1, \dots ,N-1)$. \\
b) $S^{\lambda }$ is connected, i.e. for any $\sigma \in S^{\l }$, there exist $i_1 , \dots ,i_r$ such that if we put $\s _{l} = \s (i_1, i_1+1) \dots (i_l, i_l +1) $ then $\s _{r}= min $,  $\s _{l} \in  S^{\lambda }$, $\s _{l} \succ \s _{l+1} \; (l=1, \dots ,r)$ . \\
c) Suppose $\s \in S^{\lambda }$, then $\s (i,i+1) \in S^{\lambda } \Leftrightarrow \l _{\s (i)} \neq  \l _{\s (i+1)}$. \\
d) If $\l _{\s (i)} > \l _{\s (i+1)}$ and $\s \in S^{\lambda }$ then $\s \succ \s (i,i+1)$.  

We define the partial ordering of the set $\{ (\lambda , \sigma )\: | \: \l \in \MN, \sigma \in S^{\lambda }\}$:
\begin{equation}
  (\lambda , \sigma ) \succ ( \tilde{\lambda } , \tilde{\sigma }) 
  \; \; \; \Leftrightarrow \; \; \; 
|\lambda |= |\tilde{\lambda } | \mbox{ and }
\left\{ 
\begin{array}{l} 
\lambda \succ \tilde{\lambda }   \\
\lambda = \tilde{\lambda } , \; \; \sigma \succ \tilde{\sigma} . 
\end{array} \right. 
\end{equation}

Then $Y_{i}^{(N)}$ act triangularly on $\cz $
 with respect to this ordering \cite{U1}:
\begin{equation}
 Y_{i}^{(N)} z^{\lambda ^{\sigma }} = \xi _{i}^{\lambda}(\sigma) 
z^{\lambda ^{\sigma }} + \mbox{``lower terms''} , 
\label{Ytri}
\end{equation}
\vspace{-.2in}
\begin{equation}
\xi _{i}^{\lambda}(\sigma) = 
p^{\lambda _{\sigma (i)}} q^{2\sigma (i) -N-1}  \; \; \; \;
(\sigma \in S^{\lambda } ). 
\end{equation}
In the above notation, we identify the ordering of  monomials  
$z^{\lambda ^{\sigma }}:= z_{1}^{\l _{\s_{(1)}}} z_{2}^{\l _{\s_{(2)}}}
\dots z_{N}^{\l _{\s_{(N)}}}$ with the ordering on the set of pairs $(\l, \s)$.

For generic $q$ and  $p$  the pair $(\lambda , \sigma )$ is uniquely determined from the ordered set
 $ (\xi _{1}^{\lambda}(\sigma) , \xi _{2}^{\lambda}(\sigma) , \dots
\xi _{N}^{\lambda}(\sigma) )$:
 \begin{equation}
 (\lambda , \sigma ) \neq  (\tilde{\lambda } , \tilde{\sigma})
 \; \; \; \Leftrightarrow \; \; \; 
( \xi _{1}^{\lambda}(\sigma) , \xi _{2}^{\lambda}(\sigma) , \dots
\xi _{N}^{\lambda}(\sigma) ) \neq 
( \xi _{1}^{\tilde{\lambda}}(\tilde{\sigma }) , 
\xi _{2}^{\tilde{\lambda }}(\tilde{\sigma }) , \dots
\xi _{N}^{\tilde{\lambda }}(\tilde{\sigma }) ).
\end{equation}
 Therefore one can simultaneously diagonalize
 the operators $Y_{i}^{(N)} \; \; (1 \leq i \leq N)$.
\begin{equation}
 Y_{i}^{(N)} \Phi_{\sigma}^{\lambda } (z)=
\xi _{i}^{\lambda}(\sigma) \Phi_{\sigma}^{\lambda } (z),\; \; \; \; 
\Phi_{\sigma}^{\lambda } (z)=z^{\lambda ^{\sigma }}
+ \mbox{``lower terms''} . 
\label{Ydia}
\end{equation}
 The Laurent polynomial $\Ph $ is known as  the non--symmetric Macdonald polynomial.

 The action of $g_{i,i+1}$ on the non--symmetric Macdonald polynomial is 
as follows \cite{U1}.
\begin{equation}
 g_{i,i+1} \Phi _{\sigma }^{\lambda } (z)=
A_{i}(\sigma ) \Phi _{\sigma }^{\lambda } (z)+   
B_{i}(\sigma ) \Phi _{\sigma (i,i+1)}^{\lambda }(z),
\end{equation}
\vspace{-.3in}
\begin{equation}
 A_{i}(\sigma ):=\frac{(q-q^{-1})x}{x-1} , \; \; \; \;
B_{i}(\sigma ):= 
\left\{ 
\begin{array}{cc}
 q^{-1} \{ x \} \; \; & (\lambda _{\sigma (i)} > \lambda _{\sigma (i+1)}); \\
 0 & (\lambda _{\sigma (i)} = \lambda _{\sigma (i+1)}); \\
 q^{-1} & (\lambda _{\sigma (i)} < \lambda _{\sigma (i+1)}), 
\end{array} \right. 
\end{equation}
\vspace{-.2in}
\begin{equation} 
\{ x \} :=\frac{(x-q^{2})(q^{2}x-1)}{(x-1)^{2}} ,\; \; \; \;
x:=\frac{\xi _{i+1}^{\lambda}(\sigma)}{\xi _{i}^{\lambda}(\sigma)}. 
\label{Ydiae}
\end{equation}

The case $p=1$ is not generic. However, from  the  results of  \cite{Knop} it follows that the coefficients of $\Ph$ have no poles at $p=1.$ Therefore  the non--symmetric Macdonald polynomials $\Ph $ are still well--defined at $p=1$ and the formulas (\ref{Ydia}--\ref{Ydiae}) are still satisfied.

In what follows we will let  $\Php $ denote the non--symmetric Macdonald polynomial at $p=1$. 
In virtue of the triangularity (\ref{Ytri}) the non--symmetric Macdonald polynomials $\Php \; ( \lambda \in \MN $, $ \sigma \in S^{\lambda })$ form a base of $\cz $. We put
\begin{equation}
E^{\lambda }= \bigoplus _{\s \in S^{\lambda }} \cplx \Php .
\end{equation}
Then  $\cz = \oplus _{\lambda } E^{\lambda }$.
In Section \ref{dc} we will use the following lemma.
\begin{lemma} \label{p1nsMac}
Let  $e_{-k}= \sum _{1 \leq n_1 < \dots < n_k \leq N} z_{n_1}^{-1} \dots z_{n_k}^{-1}.$ Suppose that  $\l \in \MN $ satisfies  $\l _i - \l _{i+1} = 0$ or $1.$ Then we have
\begin{equation}
e_{-k} \tilde{\Phi }_{\varsigma } ^{\l } (z) = \tilde{\Phi }_{\varsigma } ^{\tilde{\l }} (z).
\end{equation}
Here $\tilde{\l } = (\l _{1} ,\dots ,\l _{N-k} ,\dots ,\l_{N-k+1}-1 ,\dots ,\l _{N}-1)$ and $\varsigma (\in S^{\l }, S^{\tilde{\l }})$ is the minimal element of $S^{\l }$.
\end{lemma}
\begin{pf}
By the triangularity of the non--symmetric Macdonald polynomial (\ref{Ytri}), we have
\begin{eqnarray}
& & e_{-k} \tilde{\Phi }_{\varsigma } ^{\l } (z)  = e_{-k} ( z^{\l ^{\varsigma }} + \sum_{\mu \prec \l , \; \s \in S^{\mu }} c _{\mu , \s } z^{\mu ^{\s }})  
\label{eP} \\
& & = z^{\tilde{\l } ^{\varsigma }} + \sum_{(\mu , \s ) \prec (\tilde{\l }, \varsigma )} c' _{\mu , \s } z^{\mu ^{\s }} 
= \tilde{\Phi }^{\tilde{\l }} _{\varsigma } (z) + \sum_{(\mu , \s ) \prec (\tilde{\l }, \varsigma )} c'' _{\mu , \s } \tilde{\Phi } ^{\mu } _{\sigma } (z).
\nonumber
\end{eqnarray}
At $p=1$, the operators $Y_{i}^{(N)}$ commute with symmetric Laurent polynomials considered as multiplication operators on $\cz. $ Hence  we have 
\begin{equation}
Y_{i}^{(N)} e_{-k} \tilde{\Phi }_{\varsigma } ^{\l } (z) = e_{-k}Y_{i}^{(N)} \tilde{\Phi }_{\varsigma } ^{\l } (z) = q^{2 \varsigma (i)-N-1}e_{-k} \tilde{\Phi }_{\varsigma } ^{\l } (z) .
\label{jev}
\end{equation}
The ordered set of eigenvalues $\{q^{2 \sigma (i)-N-1} \}_{i=1}^N$ determines the element $\sigma \in S^{\tilde{\l}}$ uniquely. Hence (\ref{eP}) and (\ref{jev}) lead to 
\begin{equation}
e_{-k} \tilde{\Phi }_{\varsigma } ^{\l } (z) =\tilde{\Phi }^{\tilde{\l }} _{\varsigma } (z) + \sum_{\mu \prec \tilde{\l } } c'' _{\mu  } \tilde{\Phi } ^{\mu } _{\varsigma } (z).
\label{eP2}
\end{equation}
Now let us  consider any  $\mu $ which appears in the sum (\ref{eP2}).
If there exists $i (< N-k )$ such that $\mu _i < \l _i$ then for $j >i$ we necessarily  have $\mu _j < \l _j $ because of the assumption $\l _i - \l _{i+1}=0 $ or $1$ and the fact that $\l _i < \l _j $ implies  $\mu _i < \mu _j $, which follows since  $\varsigma \in S^{\mu }$.
But  $\mu _j < \l _j \; (j >i )$ leads to $|\tilde{\l }| > |\mu |$ which is in contradiction with  $ \mu \prec \tilde{\l }.$ 
Thus we have $\mu _i \geq \l _i \; (i <N-k).$ If $\mu _i > \l _i $ for some $i$, then necessarily $ \mu \not \preceq \tilde{\l }$ which is again a contradiction. Hence  we get $\mu _i = \l _i \; (i <N-k)$. 
By the condition $ \mu \prec \tilde{\l }$, we have $\mu _{N-k} = \tilde{\l } _{N-k}$. Because of the assumption on $\l $ and the fact that $\l _i < \l _j $ implies  $\mu _i < \mu _j $, we get $\mu _j \leq \tilde{\l } _j \; (j >N-k)$.
Combining this  with $|\tilde{\l }| = |\mu |$ we conclude that  $\mu _j = \tilde{\l } _j$ $\forall j$. That is  the second term in  the r.h.s. of (\ref{eP2}) is zero. 
\end{pf}

\subsection{The decomposition} \label{dc}
Let us consider the quotient space  $F_{M} / H_{-}' F_{M}$ and for each $k \geq 0$ its subspace $F^k_{M} / (H_{-}' F_{M} \cap F^k_{M}).$  

It is straightforward to establish the necessary and sufficient condition for the  vector $\omega = \sum _{\lambda } \sum _{\sigma \in S^{\lambda }} \Php \otimes \psi ^{\lambda } _{\sigma }$ $ ( \in \cz \otimes (\otimes ^{N} V))$ to be equivalent to $0$ in the quotient space $\wedge ^N V(z)$. The result is 
\begin{equation}
\forall \lambda
\left\{
\begin{array}{ll}
\psi ^{\lambda } _{\sigma (i,i+1)} = - \check{R}_{i,i+1} (q^{2(\s (i)-\s (i+1))})\psi ^{\lambda } _{\sigma }, &
\forall \sigma \mbox{ s.t. } \lambda _{\sigma (i)} > \lambda _{\sigma (i+1)} ; \\
(q^{-2}S_{i,i+1}^{-1} - S_{i,i+1} ) \psi ^{\lambda } _{\sigma }=0  , & 
\forall \sigma \mbox{ s.t. } \lambda _{\sigma (i)} = \lambda _{\sigma (i+1)} ,
\end{array}
\right.
\label{phirel1}
\end{equation}
where $\check{R}_{i,i+1} (x)$ is defined in (\ref{Rcheck}). In view of  the properties of the set $S^{\lambda }$, in the space $\wedge ^{N} V(z)$ we have
\begin{equation}
\Php \otimes \psi  _{\sigma } 
\label{Pp1Pmin}
\end{equation}
\vspace{-.2in}
\begin{equation}
\sim  \Phm  \otimes \check{R}_{i_r,i_r +1} (q^{2(\s _{r}(i_r +1) - \s _{r}(i_r))}) \dots \check{R}_{i_1,i_1 +1} (q^{2(\s _{1}(i_1 +1) - \s _{1}(i_1))}) \psi _{\sigma }.
\nonumber
\end{equation}
Here we used the notations of Section \ref{nsMac}. 
\vspace{.2in}

By the triangularity of the non--symmetric Macdonald polynomial (\ref{Ytri}) and the relation (\ref{Pp1Pmin}), we get
\begin{equation}
V_{M}^{s+nk,k} = \bigoplus _{\l } (E^{\l } \otimes (\otimes ^{N} V)) / \Omega \cap (E^{\l } \otimes (\otimes ^{N} V)) 
\label{red0}
\end{equation}
\vspace{-.2in}
\begin{equation}
= \bigoplus _{\l }(\Phm  \otimes (\otimes ^{N} V)) / \Omega \cap (\Phm \otimes (\otimes ^{N} V)),
\nonumber
\end{equation}
where the summation  is over $\lambda \in \MN ^n$, such that  $\l_{1} \leq m_{s+nk}^0$, $| \mm ^0 - \l ^{min} | =k$.

\begin{prop} \label{red1}
Define the set $\tilde{\cal M}_{s+nk}^{n,k}$ as  
\begin{equation}
\tilde{\cal M}_{s+nk}^{n,k} = \{ \l \in \tilde{\cal M}_{s+nk}^n | \:  \l_{1} \leq m_{s+nk}^0, \; | \mm ^0 - \l ^{min} | =k \mbox{ and }\l _i -\l _{i+1} =0 \mbox{ or } 1 \} .
\end{equation}
Every vector from the linear space  $F^k_{M} / (H_{-}' F_{M} \cap F^k_{M}) $ can be expressed as a linear combination of vectors of the form $\wedge (\Phm  \otimes \psi ^{\l }) |M-s-nk \rangle $, where  $\lambda \in \tilde{\cal M}_{s+nk}^{n,k}$ and $\psi ^{\l } \in  \otimes ^{N} V.$ 
\end{prop}
\begin{pf}
By the equation (\ref{red0}), it is sufficient to show that $\wedge( \Phm \otimes \psi ^{\l}) |M-s-nk \rangle $ ($\l \in \tilde{\cal M}_{s+nk}$, $\l_{1} \leq m_{s+nk}^0, \;  | \mm ^0 - \l ^{min} | =k $, $\psi ^{\l } \in  \otimes ^{N} V$ ) is equivalent to 0 in the space $F^k_{M} / (H_{-}' F_{M} \cap F^k_{M})$ unless $\l _i -\l _{i+1} =0 $ or $1$ for all $ i=1,\dots,N-1$  . We will prove this by  induction with respect to the ordering of the set $\tilde{\cal M}_{s+nk}$. (Note  that if $\l $ is not $n$--strict then $\wedge (\Phm \otimes \psi ^{\l} )= 0 $.)

Since $\l _i - \l _{i+1} \neq 0,1 $ implies that  $(\l_1, \dots ,\l_i -1, \l_{i+1}+1, \dots )$ is lower with respect to the ordering of $\tilde{\cal M}_{s+nk}$,
the minimal element satisfies the condition of  Proposition \ref{red1}.

Fix $\l $ and assume that the proposition is  proved for all $\mu $ such that  $\mu \prec \l $. Define $\tilde {\l } \in \tilde{\cal M}_{s+nk}$ as follows:
\begin{eqnarray}
& \l _i = \l _{i+1} \; \Leftrightarrow \; \tilde {\l }_i  = \tilde {\l }_{i+1}, & \\
& \l _1 = \tilde {\l }_1 ,& \\
& \tilde {\l }_i  \neq \tilde {\l }_{i+1} \; \Rightarrow \; \tilde {\l }_i  = \tilde {\l }_{i+1} +1 .& 
\end{eqnarray}
For each positive integer $l$, define $n_l = \# \{ i \: | \:  \tilde {\l }_i - \l _i =l \}$. If $n_l =0$ for all $l( \geq 1)$ then either $\wedge(\Phm \otimes \psi ^{\l} )| M-s-nk \rangle$ itself satisfies the condition of the Proposition \ref{red1}, or  else the element $\prod _{l \geq 1} B_{-l}^{n_l} \cdot \wedge(\tilde{\Phi }^{\tilde{\l}}_{min} (z)\otimes \psi _{\l }) | M-s-nk \rangle $ is in $H_{-}' F_M \cap F_M^k$.
Expanding the last element, in the space $F^k_{M} / (H_{-}' F_{M} \cap F^k_{M}) $ we get
\begin{equation}
\wedge (\Phm \otimes \psi ^{\l } ) | M-s-nk \rangle+ \sum_{\mu \prec \l , \; \s \in S^{\mu }} \wedge (\tilde{\Phi }^{\mu }_{\s } (z)\otimes \psi ^{\mu }_{\s })| M-s-nk \rangle \sim 0,
\end{equation}
for some $\psi ^{\mu }_{\s }$. By (\ref{Pp1Pmin}) and the induction assumption the proposition is proven.
\end{pf}

\begin{prop} \label{MNnk}
For each $\lambda \in \tilde{\cal M}_{s+nk}^{n,k}$, define $J$ and $r_j$ such that 
$\l _1 = \dots = \l _{r_J} >  \l _{r_J +1} = \dots =\l _{r_J +r_{J-1}} > \dots \geq \l _{r_1 + \dots +r_{J}(=N)}$, then in the space $F^k_{M} / (H_{-}' F_{M} \cap F^k_{M}) $, for each $\psi \in \otimes ^N V$ we have  
\begin{equation}
\wedge (\Phm  \otimes R_{i,i+1}(q^{2} )\psi ) |M-s-nk \rangle
\sim 0 , \; \; \; \; ( \l _{min (i)} = \l _{min (i+1)}) 
\end{equation}
\vspace{-.2in}
\begin{equation}
\wedge (\Phm  \otimes \prod_{1 \leq a \leq r_j \atop{ 0 \leq b \leq r_{j+1}-1}} R_{l_j +a, l_j+r_j+r_{j+1}-b}(q^{-2(a+b)})  \psi ) |M-s-nk \rangle \sim 0 ,
\label{jrel}
\end{equation}
where $(a,b)$ is on the right to $(a',b')$ in the product if $a<a'$ or ($a=a'$ and $b<b'$), $l_j = \sum_{i=1}^{j-1} r_j$ and $l_0 =0$.
\end{prop}
\begin{pf}
The first relation follows from  (\ref{phirel1}) and the identity 
\begin{equation}
\mbox{Im} (q^{2} S_{i,i+1}^{-1} - S_{i,i+1} ) = \mbox{Ker} (q^{-2} S_{i,i+1}^{-1} - S_{i,i+1} ) .
\end{equation}
Consider the second relation. We define  
\begin{equation}
\bar{\l } = ( \l_1 , \dots , \l _{r_{j+1}+ \dots +r_{J}}, \l _{1+r_{j+1}+ \dots +r_{J}}+1, \dots ,\l_{N}+1).
\end{equation}
By Lemma \ref{p1nsMac}, the definition of the space $F^k_{M} / (H_{-}' F_{M} \cap F^k_{M}) $ and the relation (6.51) in \cite{STU}, we get
\begin{equation}
f( B_{-1}, \dots , B_{-(r_1 + \dots +r_j)} )\cdot \wedge (\tilde{\Phi }^{\bar{\l }}_{min} (z)\otimes \psi ) |M-s-nk \rangle 
\end{equation}
\vspace{-.2in}
\begin{equation}
= ( f( B_{-1}, \dots , B_{-(r_1 + \dots +r_j)} ) \cdot \wedge (\tilde{\Phi }^{\bar{\l }}_{min} (z) \otimes \psi ))|M-s-nk \rangle 
\nonumber
\end{equation}
\vspace{-.2in}
\begin{equation}
=\wedge (\tilde {\Phi }^{\l }_{\varsigma } (z) \otimes \psi )| M-s-nk \rangle \sim 0.
\nonumber
\end{equation}
Here $f( x_1 , \dots , x_l)$ is a polynomial such that 
\begin{equation}
f (\sum _{i=1}^{N} z_i , \sum _{i=1}^{N} z_i^2 ,\dots ,\sum _{i=1}^{N} z_i^{l} ) = \sum _{i_1 < \dots < i_l } z_{i_1} \dots z_{i_l},
\end{equation}
 and $\varsigma \in S^{\l} $ is the minimal element of $S^{\bar{\l}}$.

If we apply the formula (\ref{Pp1Pmin}), we get
\begin{equation}
\wedge (\tilde {\Phi }^{\l }_{min} (z) \otimes  \prod_{0 \leq a \leq r_j-1 \atop{ 0 \leq b \leq r_{j+1}-1}} \check{R}_{l_j +r_{j+1}-b+a, l_j +r_{j+1}-b+a+1}(q^{-2(a+b+1)}) \psi )| M-s-nk \rangle \sim 0,
\end{equation}
where $(a,b)$ on the right to $(a',b')$ in the product if $a<a'$ or ($a=a'$ and $b<b'$).
Finally, taking into account the relation
\begin{equation}
 \prod_{0 \leq a \leq r_j-1 \atop{ 0 \leq b \leq r_{j+1}-1}} \check{R}_{l_j +r_{j+1}-b+a, l_j +r_{j+1}-b+a+1}(q^{-2(a+b+1)}) \cdot \prod_{1 \leq a \leq r_j \atop{ 0 \leq b \leq r_{j+1}-1}} P_{l_j +a, l_j+r_j+r_{j+1}-b}
\end{equation}
\vspace{-.2in}
\begin{equation}    
 = \prod_{1 \leq a \leq r_j \atop{ 0 \leq b \leq r_{j+1}-1}} R_{l_j +a, l_j+r_j+r_{j+1}-b}(q^{-2(a+b)}),
\nonumber
\end{equation}
we obtain  (\ref{jrel}).
\end{pf}

With notations of  Proposition \ref{MNnk}, for each  $\l \in \tilde{\cal M}_{s+nk}^{n,k}$ define the linear subspace of $\otimes ^N V$ 
\begin{equation}
V^{\l } = \sum _{\l _{min(i)} = \l _{min(i+1)}} \mbox{Im} R_{i,i+1}(q^2) \; + \sum_{j=1}^{J-1} \mbox{Im} (\prod_{1 \leq a \leq r_j \atop{ 1 \leq b \leq r_{j+1}-1}} R_{l_j +a, l_j+r_j+r_{j+1}-b}(q^{-2(a+b)})).
\end{equation}

By Proposition \ref{MNnk}, we get
\begin{prop} \label{Psurj}
Consider the map 
\begin{eqnarray}
\psi_k : \bigoplus _{\l \in  \tilde{\cal M}_{s+nk}^{n,k}} \Phm \otimes (\otimes ^{s+nk} V  / V^{\l }) & \rightarrow & F_M^k / ( H_{-}'F_M \cap F_M^k ) \\
\label{psik}
v & \mapsto & v \wedge |M-s-nk \rangle . \nonumber
\end{eqnarray}
Define the action of $U'_q (\sll)$ on the l.h.s. of {\em (\ref{psik})} by {\em (\ref{e: Efin}-\ref{e: K0})}.

Then  the map $\psi_k $ is well-defined, surjective and is a  $U'_q (\sll)$--intertwiner.
\end{prop}

\begin{prop} \label{lteq}
Let $\l \in \tilde{\cal M}_{s+nk(=N)}^{n,k}.$ For any such  $\l$ we define  $J$, $r_j$ and $l_j$ in the same way  as in Proposition \ref{MNnk}.
 Let $\theta $ be the border strip characterized by $\langle r_J, r_{J-1} , \dots , r_1 \rangle $. We have
\begin{equation}
\otimes ^N V / V^{\l } = \otimes ^N V / Ker \bar{R}_{\theta }.
\end{equation}
\end{prop}
\begin{pf}
First we will show that $\otimes ^N V / V^{\l } \supset \otimes ^N V / Ker \bar{R}_{\theta }$.
Applying repeatedly the Yang--Baxter equation $R_{a,b} (x) R_{a,c}(xy) R_{b,c}(y)= R_{b,c}(y)R_{a,c}(xy)R_{a,b} (x)$, we can move some special elements to the right in the product 
\begin{equation}
\bar{R}_{\theta } = \cdots R_{i+1,i}(q^{-2}) 
 = \cdots \prod_{ 1 \leq a \leq r_{j+1} \atop{ 0 \leq b \leq r_j -1}} R_{ l_j+r_j+a, l_j +r_j-b }(q^{2(a+b)}),
\end{equation}
where $(a,b)$ on the right to $(a',b')$ in the product if $a<a'$ or ($a=a'$ and $b<b'$), and $\l _{min(i)} = \l _{min(i+1)}$.
By the  formula
\begin{equation}
R_{b,a} (x) R_{a,b} (x^{-1}) = \frac{(x-q^2)(x^{-1}q^{-2}-1)}{(x-1)(x^{-1}-1)} id,
\end{equation}
we get $\otimes ^N V / V^{\l } \supset \otimes ^N V / Ker \bar{R}_{\theta }$.

Next we will show that $\otimes ^N V / V^{\l } \subset \otimes ^N V / Ker \bar{R}_{\theta }$. 
We show that the vectors $\otimes  _{i=1}^N v_{e_i}$ such that $e_i <e_{i+1} $ if $\l _{min{(i)
}} = \l _{min{(i+1)}}$, $e_i \geq e_{i+1} $ if $\l _{min{(i)}} \neq \l _{min{(i+1)}}$  span the space $\otimes ^N V / V^{\l }$.

Using the relations $R_{i',i'+1}(q^{2})( \otimes ^N v_{e_i} ) \sim 0$ $(\l _{min(i')}= \l _{min(i'+1)})$ we find that the set of vectors $\{ \otimes _{i=1} ^N v_{e_i}\: |\: \text{if $\l _{min{(i)}} = \l _{min{(i+1)}}$ then $e_i <e_{i+1}$ }\} $ spans the space  $\otimes ^N V / V^{\l }$. For each  vector of the form $ \otimes  _{i=1}^N v_{e_i}$ we define $\tilde{N} (\otimes  _{i=1}^N v_{e_i} ) = \# \{ (i,j)| \; i<j,\; e_i \geq e_j $ and $\l _{min{(i)}} \neq \l _{min{(j)}} \} $.
 Consider the vector $\otimes  _{i=1}^N v_{e_i} $ such that if $\l _{min{(i)}} = \l _{min{(i+1)}}$ then $e_i <e_{i+1} $. 
Assume that there is an $i$ such that $e_i < e_{i+1} $ (if $\l _{min{(i)}} \neq \l _{min{(i+1)}}$), then we get the relation 
\begin{equation}
\prod_{1 \leq a \leq r_j \atop{ 1 \leq b \leq r_{j+1}-1}} R_{l_j +a, l_j+r_j+r_{j+1}-b} ( q^{-2(a+b)} )   (  \otimes ^N v_{e_i} )  \sim 0,
\label{re0}
\end{equation}
for all possible $j$. By (\ref{re0}), the vector $\otimes ^N v_{e_i} $ is equivalent to a linear combination  of the vectors $\otimes ^N v_{e_{i'}} $ such that $\tilde{N} (\otimes ^N v_{e_{i'}} ) < \tilde{N} (\otimes ^N v_{e_i} )$.
Because $\otimes ^N v_{e_{i'}} $ is invariant by the relations $R_{i',i'+1}(q^{2})( \otimes ^N v_{e_i} )$ ($ l_j +1 \leq i' \leq l_{j+1}-1, \;  l_{j+1} +1 \leq i' \leq l_{j+2}-1$), if we use these relations we get that the vector  $\otimes ^N v_{e_i} $ ($e_i <e_{i+1} $ if $\l _{min{(i)}} = \l _{min{(i+1)}}$) is expressed by the sum of $\otimes ^N v_{e_{i'}} $ such that $\tilde{N} (\otimes ^N v_{e_{i'}} ) < \tilde{N} (\otimes ^N v_{e_i} )$ and $e_{i'} < e_{i'+1} $ (if $\l _{min{(i')}} \neq \l _{min{(i'+1)}}$).

By the induction on $\tilde{N}$ we find  that the vectors $\otimes ^N v_{e_i}$ ($e_i <e_{i+1} $ if $\l _{min{(i)}} = \l _{min{(i+1)}}$, $e_i \geq e_{i+1} $ if $\l _{min{(i)}} \neq \l _{min{(i+1)}}$) span the space $\otimes ^N V / V^{\l }$.

The number of these vectors is equal to the number of semi--standard tableaux  of $\theta $. By Corollary \ref{dimbs}, we get $\otimes ^N V / V^{\l } \subset \otimes ^N V / Ker \bar{R}_{\theta }$.
\end{pf}
\vspace{.2in}

In what follows we  identify the  border strips $\langle m_1, \dots , m_r \rangle $ and $\langle m_1, \dots , m_r ,n, \dots ,n \rangle $, and identify $(\l_1 , \dots ,\l _{s+nk} )$ and $(\overbrace{\l_1 +1, \dots ,\l_1 +1}^{n} ,\l_1 , \dots ,\l _{s+nk})$ which are elements of $\coprod _k \tilde{\MC }^{n,k}_{s+nk}$. Proposition \ref{lteq} gives a one to one correspondence of  $\coprod _k \tilde{\MC }^{n,k}_{s+nk}$ and the set of all skew Young diagrams of the border strip type  $\langle m_1, \dots , m_r \rangle $  which satisfy $m_i \leq n $ for all $i$ and $\sum_{i=1}^{r} m_i$ $\equiv$ $s$ mod $n$. 
On this correspondence the degree of the semi--infinite wedge $ \wedge( \Phm \otimes \psi ) | M-s-nk \rangle$ is equal to $\frac{1-n}{24} - \frac{k(n-k)}{2n} + \frac{1}{2n} | \theta |(n- | \theta |)+t(\theta ) $, where $\theta $ is the  border strip which corresponds to $\l $ and $t(\theta ) = \sum_{i=1}^{r-1} (r-i) m_i$. 

We define $ch (F_M / H_- ' F_M)= \sum _{\mu ,i } dim(V_{\mu ,i}) e^{\mu } q^i$, where $V_{\mu ,i}$ is the subspace of $F_M / H_- 'F_M$ of the degree (\ref{degifw}) $i$ and of the  $U_q(\sln )$--weight $\mu $.
We put $ \sum _{\mu ,i } a_{\mu ,i} e^{\mu } q^i \leq \sum _{\mu ,i } b_{\mu ,i} e^{\mu } q^i$ iff $a_{\mu ,i} \leq b_{\mu ,i}$ for all $\mu $ and $i$.
By Proposition \ref{Psurj} we have
\begin{equation}
ch (F_M / H_- ' F_M) \leq q^{\frac{1-n}{24}- \frac{k(n-k)}{2n}} \sum_{\theta \in BS \atop {| \theta | \equiv k \; mod \; n}} q ^{\frac{1}{2n} | \theta | (n- | \theta |)+t(\theta )} s_\theta (z).
\label{chFM}
\end{equation}
$F_M / H_- ' F_M$ with $U_1$--action is isomorphic to $V(\Lambda _k) $, this isomorphism is degree preserving with respect to the degree (\ref{degifw}) on $F_M / H_- ' F_M$ and the homogeneous degree on $V(\Lambda _k) $, and the character formula of $V(\Lambda _k)$ is given in Proposition \ref{chara}. Hence the inequality of (\ref{chFM}) is, in fact, an  equality and, therefore,  the map (\ref{psik}) must is  bijective. Thus have the following theorem  
\begin{thm}
We have the isomorphism of  $U'_q (\sll )$--modules:
\begin{equation}
F_M / H_- ' F_M \simeq \bigoplus _{\theta } V_{\theta }, \label{deco}
\end{equation}
where the sum is over all border strips $\langle m_1, \dots m_r \rangle $, $(m_i \leq n $, $m_r < n$ and $N\equiv M $ mod $n)$, the space $V_{\theta }$ and the level--$0$ $U'_q (\sll )$--action is defined by $( \pi ^{(N)} _{a_1, \dots , a_{N}} ,\; R_{\theta } \cdot \otimes ^{N} V )$ where $N= \sum_{i=1}^{r} m_i$ and $a_{l+\sum_{i=j}^{r} m_i} = 2(l-1+\sum_{i=1}^{j-1} m_i)$.

\end{thm}
 
\subsection{$\tl $ case}

In this section we will discuss the $\tl $ case in a somewhat more detail.

Let $W_n$ be the $(n+1)$--dimensional irreducible module of $U_q(\tl )$, and $W_n(b)$ be the evaluation module with the parameter $b$ whose $\Uslt $--module structure is given by
\begin{equation}
E_0 = q^b F_1, \; \; \; F_0 = q^{-b}E_1, \; \; \; K_{0}= K_{1}^{-1}.
\end{equation}

It is known that every finite dimensional irreducible $\Uslt $--module is isomorphic to $\otimes _{\mu } W_{n_{\mu }} (b_{\mu})$ for some $n_{\mu }$ and $ b_{\mu}$. We will represent the $\Uslt $--module described by a skew Young diagram as the tensor product of the form  $\otimes _{\mu } W_{n_{\mu }} (b_{\mu})$.

\begin{prop} \label{sltP}
Let $\theta $ be the skew Young diagram of border strip $\langle m_1, \dots, m_r \rangle $ such that $m_i = 1$ or $2$, and $N= \sum_{i=1}^{r} m_i$, $a_{l+\sum_{i=j}^{r} m_i} = 2(l-1+\sum_{i=1}^{j-1} m_i)$.
We put $I= \{ i \: | \: m_i =1 $ and $ m_{i-1} =2\} = \{ l_1 , \; l_2 , \dots , l_{r'} \}$ $(m_0=2, l_i< l_{i+1})$ and let $n_i$ be the integer such that $m_{l_i} = \dots = m_{l_i +n_i -1} =1, \; m_{l_i +n_i} =2$ and $b_i= 2l_i+n_i -3$.

The $\Uslt $--module $(  \pi ^{(N)} _{a_1, \dots , a_{N}} , \: R_\theta \cdot \otimes^N \cplx ^2  )$ is isomorphic to $W_{n_1}(b_1) \otimes W_{n_2}(b_2) \otimes \dots \otimes W_{n_{r'}}(b_{r'})$.
\end{prop}

\begin{pf}
By Proposition \ref{skewrep}, we find that  $( \pi ^{(N)} _{a_1, \dots , a_{N}},\: R_\theta  \cdot \otimes^N \cplx ^2 )$ is isomorphic to $( \pi ^{(N)} _{a_N, \dots , a_1} , \: \otimes^N \cplx ^2 / $Ker$\bar{R}_{\theta })$. 

As in  the proof of Proposition \ref{lteq}, we get Im$\check{R}_{i,i+1} (q^{2}) \subset $Ker$\bar{R}_{\theta }$ if $a_{i+1} = a_i -2$ and Im$ \check{R}_{i,i+1} (q^{-2}) \subset $Ker$\bar{R}_{\theta }$ if $a_{i+1} = a_i +2$. 

We can directly confirm that the $\Uslt $--module $( \pi _{a,a-2}, \: \cplx ^2 \otimes \cplx ^2 / $Im$ \check{R}_{1,2} (q^{2}) )$ is $1$--dimensional and the module $( \pi _{a,a+2,\dots ,a+2(l-1)}, \: \otimes ^l \cplx ^2 / \sum_{i=1}^{l-1} $Im$ \check{R}_{i,i+1} (q^{-2}))$ is isomorphic to $W_{l}(a+l-1)$, where $(l)$ is the Young diagram of degree $l$ which has only one row.

If we put 
\begin{equation}
\tilde{V} = \sum _{ i | a_{i+1} = a_i +2 } Im \check{R}_{i,i+1} (q^{2})
+ \sum _{ i | a_{i+1} = a_i -2  } Im \check{R}_{i,i+1} (q^{-2}),
\end{equation}
then  $(\pi ^{(N)} _{a_N, \dots , a_1}, \: \otimes^N \cplx ^2 / \tilde{V}) \simeq W_{n_1}(b_1) \otimes W_{n_2}(b_2) \otimes \dots \otimes W_{n_{r'}}(b_{r'})$.

Since the dimension of $W_{n_1}(b_1) \otimes W_{n_2}(b_2) \otimes \dots \otimes W_{n_{r'}}(b_{r'})$ is equal to the dimension $\otimes^N \cplx ^2 / $Ker$\bar{R}_{\theta }$ the proof is finished.
\end{pf}
\vspace{.2in}

By Proposition \ref{sltP}, we can rewrite the decomposition (\ref{deco}) for $\tl $ case. In fact we get the same decomposition as \cite{J}. More precisely we get
\begin{prop}
If we change the coproduct of the level--0 $\Uslt $--module defined in \cite{J} to fit our coproduct (\ref{cp1}--\ref{cp3}), The level--0 $\Uslt $--module of $V(\Lambda _s)$ $(s \equiv M $ mod $2, \; s=0$ or $1)$ defined in \cite{J} is isomorphic to the level--0 $\Uslt $--module of $V(\Lambda _s)$ defined in this paper, and the degree is preserved under this isomorphism.
\end{prop}

{\bf Acknowledgment}
I am grateful to D. Uglov with whom we collaborated on the papers \cite{TU,STU}. I also thank M. Jimbo, M. Kashiwara, T. Miwa, M. Nazarov, Y. Saito, V. Tarasov and colleagues at  R.I.M.S. for discussions and support.

\end{document}